\documentclass[12]{article}
\usepackage{enumerate}
\usepackage{amssymb}
\usepackage[all]{xy}

\def\buildrel#1_#2^#3{\mathrel{\mathop{\kern 0pt#1}\limits_{#2}^{#3}}}

\newcommand{\Pf}{{\em Proof}. }

\newcommand{\SPAN}{\mbox{$\mathtt{span}$}}

\newcommand{\R}{\mathbb R}

\newcommand{\g}{{\mathfrak{g}}{}}

\newcommand{\ddto}{{\frac{d}{dt}|_{0}}{}}

\newcommand{\h}{{\mathfrak{h}}{}}

\newcommand{\CO}{{\cal O}{}}

\def\cref#1{Corollary~\ref{#1}}

\newtheorem{thm}{Theorem}[section]
\newtheorem{prop}[thm]{Proposition} \newtheorem{lem}[thm]{Lemma}
\newtheorem{cor}[thm]{Corollary} \newtheorem{dfn}[thm]{Definition}
 \newtheorem{rmk}[thm]{Remark}
\newtheorem{ex}[thm]{Example} 

\def\pf{\noindent{\bf Proof.}\ }

\def\complex{{\mathbb C}}

\def\reals{{\mathbb R}}

\def\integers{{\mathbb Z}}
\def\cala{{\mathcal A}}

\def\calg{{\mathcal G}}
\def\calh{{\mathcal H}}

\def\call{{\mathcal L}}
\def\calm{{\mathcal M}}

\def\calo{{\mathcal O}}

\def\calu{{\mathcal U}}

\def\frakg{\mathfrak{g}}

\newtheorem{theorem}{Theorem}

\newtheorem{definition}[theorem]{Definition}

\def\g{\gamma}

\def\part{\partial}

\def\ts{\times}

\def\text{\hbox}

\def\build#1_#2^#3{\mathrel{
\mathop{\kern 0pt#1}\limits_{#2}^{#3}}}

\def\@nbibitem#1{\noindent \hangindent=2pc \hangafter=1
\refstepcounter{enumi}\hbox to 2pc{\arabic{enumi}.\hfil}%
\immediate\write\@auxout{\string\bibcite{#1}{\arabic{enumi}}}}
\def\numbibliography{%
\section*{REFERENCES}%
\bgroup\footnotesize
\setcounter{enumi}{0}%
\def\newblock{\hskip .11em plus.33em minus.07em}%
\let\bibitem\@nbibitem}
\def\endnumbibliography{\par\egroup}
\makeatother

\def\g{\gamma}

\def\build#1_#2^#3{\mathrel{
\mathop{\kern 0pt#1}\limits_{#2}^{#3}}}
\title{Rankin-Cohen brackets and formal quantization}
\author{Pierre Bieliavsky, Xiang Tang and Yijun
Yao}
\begin{document}
\maketitle \begin{abstract}In this paper, we use the theory of
deformation quantization to understand Connes' and Moscovici's
results \cite{cm:deformation}. We use Fedosov's method of
deformation quantization of symplectic manifolds to reconstruct
Zagier's deformation \cite{z:deformation} of modular forms, and
relate this deformation to the Weyl-Moyal product. We also show
that the projective structure introduced by Connes and Moscovici
is equivalent to the existence of certain geometric data in the
case of foliation groupoids. Using the methods developed by the
second author \cite{t1:def-gpd}, we reconstruct a universal
deformation formula of the Hopf algebra $\calh_1$ associated to
codimension one foliations. In the end, we prove that the first
Rankin-Cohen bracket $RC_1$ defines a noncommutative Poisson
structure for an arbitrary $\calh_1$ action.\end{abstract}
\section{Introduction}
In the study of transversal index theory, Connes and Moscovici
introduced a Hopf algebra, $\calh_1$, which governs the local
symmetry in calculating the index of a transversal elliptic
operator. Interestingly, Connes and Moscovici
\cite{cm:modular-heck} discovered an action of $\calh_1$ on the
modular Hecke algebras.

Inspired by this action, Connes and Moscovici found many
similarities between the theory of codimension one foliations
and the theory of modular forms. For example, they showed that the
Hopf cyclic version of the Godbillon-Vey cocycle gives rise to a
1-cocycle on $PSL(2, \mathbb{Q})$ with values in Eisenstein series
of weight 2, and that the Schwarzian 1-cocycle corresponds to an
inner derivation implemented by a level 1 Eisenstein series of
weight 4. In particular, inspired by Zagier's \cite{z:deformation}
Rankin-Cohen deformation  on modular forms, Connes and Moscovici
\cite{cm:deformation} constructed  a universal deformation formula
for an action of $\calh_1$ with a projective structure. In this
paper, we aim to reconstruct this deformation formula using
noncommutative Poisson geometry as developed by the second author
\cite{t1:def-gpd} and \cite{t2:thesis}.

The origin of the Rankin-Cohen deformation is a work of Rankin.
Rankin in 1956 described all polynomials in the derivatives of
modular forms with values again in modular forms. Based on
Rankin's work, in 1977, Cohen defined a sequence of bilinear
operations on modular forms indexed by nonnegative integer $n$,
which assigns to two modular forms, $f$ of weight $k$ and $g$ of
weight $l$, a modular form of weight $k+l+2n$. Their results
showed that for any given integer $n\geq 0$, there is
essentially(up to a constant) only one bilinear operator
mapping\footnote{$\calm_p$ is the space of modular forms of weight
$p$.} $\calm_p\otimes \calm_q$ to $\calm_{p+q+2n}$ $\forall p,
q\in \integers_{\geq0}$. They are later called Rankin-Cohen
brackets and usually denoted by $RC_n$. These operators were
further studied and played an important role in the theory of
modular forms. Zagier \cite{z:deformation} observed that the sum
of Rankin-Cohen brackets defines an associative product on the
algebra $\calm:=\sum_{l\geq0} \calm_l$. Zagier's proof of the
associativity of this product, which involves infinitely many
equalities, was rather combinatoric. Cohen, Manin, and Zagier
\cite{cmz:modular} explained this deformation using the theory of
automorphic pseudo differential operators. The calculation still
involves many interesting and complicated combinatoric identities.
In this paper, we will first reconstruct Zagier's Rankin-Cohen
deformation using the methods of deformation quantization of
symplectic manifolds developed by Fedosov \cite{fe:book}. In
particular, we will show that this deformation is isomorphic to
the  standard Moyal product. The calculation involved in our
construction is easier and more transparent than those
\cite{cmz:modular} and \cite{z:deformation}.

To reconstruct Connes-Moscovici's Rankin-Cohen deformation for
$\calh_1$ action, we need to first understand the projective
structure introduced by Connes and Moscovic \cite{cm:deformation}.
The notion of a projective structure of $\calh_1$ is a
generalization of the projective structure on an elliptic curve (
see \cite{cmz:modular}). Our idea to understand this structure is
to look at the defining action of $\calh_1$ on a groupoid algebra
associated to a codimension one foliation. In this case, we
discovered that the existence of a projective structure is
equivalent to the existence of a certain type of invariant
symplectic connection. This geometric explanation provides a
natural connection to the results in Tang \cite{t1:def-gpd}, where
he studied the deformation quantization of a groupoid algebra. The
existence of an invariant symplectic connection is a sufficient
condition for the existence of a deformation quantization of a
groupoid algebra. Therefore, in the case of a codimension one
foliation, Tang's construction \cite{t1:def-gpd} implies that with
a projective structure, one can construct a deformation
quantization (a star product) of the corresponding foliation
groupoid algebra. Furthermore, our calculation in Section 5
exhibits that when the symplectic connection is flat, the star
product on the groupoid algebra can be expressed by an element
$RC$ in $\calh_1\otimes\calh_1[[\hbar]]$. To obtain a universal
deformation for a $\calh_1$ action with a projective structure as
Connes and Moscovici \cite{cm:deformation}, we construct a fully
injective $\calh_1$ action on the union of groupoid algebras of
those foliation groupoids with a fixed type of invariant
symplectic connections. Therefore, we are able to reconstruct the
universal deformation formula on $\calh_1$ by pulling back the
star products on the groupoid algebras.

All the above deformations, including \cite{cmz:modular},
\cite{cm:deformation}, and \cite{z:deformation}, are all formal
deformation, which means that the deformation parameter $t$ is a
formal variable. It is more interesting to ask whether one can
make a deformation strict in the sense of Rieffel. This will be
studied in the next paper \cite{bty:cmz2}.

{\bf Acknowledgment}: We would like to thank Alain Connes and
Henri Moscovici for explaining the results \cite{cm:deformation}
and asking us interesting questions. Tang would like to thank Max
Karoubi and Richzard Nest for their hosts of his visit of
Institute de Henri Poincar\'e during summer 2004, where the paper
started. Yao wants to thank Don Zagier for his inspiring course
given at Coll\`ege de France.
\section{Prerequisites}

In this section, we review the materials needed for this paper.
\subsection{Codimension one foliations and the Hopf algebra}
For a constant rank foliation on $M$, we choose a complete flat
transversal $X$. We look at the oriented frame bundle $FX$ of $X$
with the lifted holonomy foliation groupoid action, which defines
an \'etale groupoid $\calg\rightrightarrows FX$. Connes and
Moscovici found a Hopf algebra $\calh_k$ acting on the smooth
groupoid algebra $C_c^{\infty}(\calg)$, where $k$ is the
codimension of the foliation. We exhibit this Hopf algebra in the
case of $k=1$.

In the case of a codimension one foliation, the complete
transversal $X$ is a flat 1-dim manifold, and $FX$ is isomorphic
to $X\times \reals^+$ by fixing a flat connection on $FX\to X$. We
introduce coordinates $x$ on the $X$ component and $y$ on the
$\reals^+$ component. Let $\Gamma$ be a pseudogroup associated to
the foliation acting on $X$. The lifted action of $\Gamma$ on $FX$
is
$$ (x,y)\mapsto (\phi(x), \phi'(x)y), \ \ \ \ \ \ \forall \phi\in \Gamma.$$

We look at the groupoid $FX\rtimes \Gamma\rightrightarrows FX$. It
is an \'etale groupoid with a natural symplectic form
$\omega=\frac{{\rm d}x\wedge{\rm d}y}{y^2}$.

On $FX$, we consider vector fields $X=y\partial_x$ and
$Y=y\partial_y$. It is easy to check that $Y$ is invariant under
the $\Gamma$ action, but $X$ is not, and has the following
commutation relation,
$$ U_{\phi}X U_{\phi}^{-1}=X-y\frac{{\phi^{-1}}''(x)}{{\phi^{-1}}'(x) }Y.$$

We introduce the following operators on $\cala$.
\begin{equation}
\label{dfn:hopf-act}
\begin{array}{ll}
X(fU_{\phi})&=X(f)U_{\phi},\\
Y(fU_{\phi})&=Y(f)U_{\phi},\\
\delta_1(fU_{\phi})&=\mu_{\phi^{-1}}fU_{\phi},\\
\delta_n(fU_{\phi})&=X^{n-1}(\mu_{\phi^{-1}})fU_{\phi},
\end{array}
\end{equation}
where
$\mu_{\phi^{-1}}(x,y)=y\frac{{\phi^{-1}}''(x)}{{\phi^{-1}}'(x)}$.

The commutation relation among the above operators are
\[
\begin{array}{ll}
[Y, X]=X,& [X, \delta_n]=\delta_{n+1},\\

[Y,\delta_n ]
 =n\delta_n,& [\delta_n, \delta_m]=0.
\end{array}
\]
The operators $X, Y, \delta_n,\ n\in \mathbb{N}$ form an infinite
dimensional Lie algebra $H_1$, and the Hopf algebra $\calh_1$ is
defined to be the universal enveloping algebra of $H_1$.

We define the following operations on $\calh_1$:
\begin{enumerate}
\item product $\cdot : \calh_1\otimes \calh_1\to \calh_1$ by the
product on $\calh_1$ as the universal enveloping algebra of $H_1$.
\item coproduct $\Delta:\calh_1\to \calh_1\otimes \calh_1$ by
$$
\begin{array}{l}
\Delta Y=Y\otimes 1+1\otimes Y,\\
\Delta \delta_1=\delta_1\otimes 1+1\otimes \delta_1,\\
\Delta X=X\otimes 1+1\otimes X+\delta_1\otimes Y,\\
\Delta\delta_n=[\Delta X, \Delta \delta_{n-1}].
\end{array}
$$
\item counit $\epsilon:\calh_1\to \complex$ by taking the value of
the identity component.
\item antipode $S: \calh_1\to \calh_1$ by
$$
S(X)=-X+\delta_1Y,\ \ S(Y)=-Y,\ \ S(\delta_1)=-\delta_1.
$$
\end{enumerate}
It is straightforward to check that $(\calh_1, \cdot, \Delta, S,
\epsilon, id)$ defines a Hopf algebra.
\subsection{Deformation quantization a la Fedosov}
Fedosov's construction of deformation quantizations of a
symplectic manifold can be formulated as follows.

Let $(M,\omega)$ be a $2n$ dimensional symplectic manifold. At
each fiber $T_xM$ of the tangent bundle, which is a symplectic
vector space, we define a Weyl algebra $W_x$ to be an associative
algebra over $\complex$ with a unit, whose elements are of the
form
\[
a(y,\hbar)=\sum\limits_{k,|\alpha|\geq 0} \hbar^k a_{k,\alpha}
y^\alpha,
\]
where $\hbar$ is a formal parameter and $y=(y^1,\dots,y^{2n})\in T_x
M$ is a tangent vector, $\alpha=(\alpha_1,\dots,\alpha_{2n})$ is a
multi-index,
$y^\alpha=(y^1)^{\alpha_1}\cdots(y^{2n})^{\alpha_{2n}}$.

The product of elements $a,b\in W_x$ is defined as follows:
\[
\begin{array}{ll}
a\circ b&=\exp(-\frac{i\hbar}{2}\omega^{ij}\frac{\partial}{\partial
y^i}\frac{\partial}{\partial z^j}) a(y,\hbar)b(z,\hbar)|_{z=y}\\
&=\sum_{k=0}^{\infty} (-\frac{i\hbar}{2})^k \frac{1}{k!}\omega^{i_1
j_1}\cdots\omega^{i_k j_k}\frac{\partial^k a}{\partial y^{i_1}\cdots
\partial y^{i_k}}\frac{\partial^k b}{\partial y^{j_1}\cdots \partial
y^{j_k}}.
\end{array}
\]

We consider the Weyl algebra bundle $W$ over $(M,\omega)$ for
which the fiber at the point $x$ is $W_x$, and denote
$C^{\infty}(W)$ to be the algebra of smooth sections of $W$ with
pointwise multiplication $\circ$. To introduce the Fedosov
connection, we look at the algebra $C^{\infty}(W\otimes
\Lambda)=\oplus_{q=0}^{2n}\Gamma^{\infty}(W\otimes \Lambda ^q)$,
where $\Lambda ^q$ is set of smooth $q-$forms.

We introduce several operations on $C^{\infty}(W\otimes \Lambda)$.
\begin{enumerate}
\item commutator, i.e. $[a, b]=a\circ b-(-1)^{deg(a)deg(b)}b\circ
a$.
\item $\delta,\ \delta^*:\ C^{\infty}(W\otimes \Lambda)\to C^{\infty}(W\otimes
\Lambda)$, i.e.
\[
\delta a=dx^k\wedge \frac{\partial a}{\partial y^k},\hspace {1cm}
\delta^* a=y^k i(\frac {\partial}{\partial x^k})a.
\]
\end{enumerate}

A \it Fedosov connection \rm on the Weyl algebra bundle $W$ is a
connection $D$ such that for any section $a\in C^\infty(W\otimes
\Lambda)$,
$$
D^2 a=\frac{i}{\hbar} [\Omega, a]=0.
$$

Fedosov in \cite{fe:book} showed that given a torsion free
symplectic connection $\nabla$ on $M$ with Christoffel
$\Gamma_{ijk}$, one can construct an abelian connection on $W$ of
the following form
$$
D=-\delta+\partial+\frac{i}{\hbar}[r,\cdot],
$$
where $\partial a:=d a+\frac{i}{\hbar}[\Gamma, a]$, with
$\Gamma=\frac{1}{2} \Gamma_{ijk} y^i y^j dx$, and $r$ is a local
1-form with values in $W$.

We look at the subalgebra $W_D\subset C^\infty(W)$ consisting of
flat sections of $D$. The main theorem that we will use is the
following:
\begin{thm}
\label{thm:fedosov} For any $a_0\in C^\infty(M)[[\hbar]]$, there
exists a unique section $a\in W_D$, which is denoted by
$\sigma^{-1}(a_0)$, such that $\sigma(a)=a_0$, where $\sigma(a)$
means the projection onto the center: $\sigma(a)=a(x,0,h)$.
\end{thm}

This implies that there is a one-to-one correspondence between
$W_D$ and $C^{\infty}(M)[[\hbar]]$. Accordingly we can define on $
C^{\infty}(M)[[\hbar]]$ an associative star product

\begin{equation}
\label{eq:star} a \star b=\sigma(\sigma^{-1}(a)\circ
\sigma^{-1}(b)).
\end{equation}

\subsection{Deformation Quantization of Groupoids}
The second named author \cite{t1:def-gpd} considered deformation
quantization of the groupoid algebra of a pseudo \'etale groupoid
and proved that one can construct star products on such groupoids.
As a special case, we have that for an \'etale groupoid with an
invariant symplectic structure and an invariant symplectic
connection on the base, the groupoid algebra can be formally
deformation quantized. In this subsection, we recall the basic
concepts and constructions from Tang \cite{t1:def-gpd}.

\begin{definition}
\label{dfn:poisson} (Block, Getzler and Xu)
  A Poisson structure on an associative  algebra $A$ is an element $[\Pi]$
  of the Hochschild cohomology group $H^2(A, A)$ such
  that the cohomology class of the Gerstenhaber bracket
  $[\Pi, \Pi]$ vanishes.
\end{definition}

\begin{definition}
  Let $(A, [\Pi])$ be a noncommutative Poisson algebra,  and $A[[\hbar]]$ the
  space of formal power series with coefficients in $A$. A
  {\it formal deformation quantization} of $(A, [\Pi])$
  (or in other words {\it star product})
  is an associative product
\[
  \star : A[[\hbar ]] \times A[[\hbar ]] \rightarrow A[[\hbar ]],
  \quad (a_1,a_2) \mapsto a_1 \star a_2 = \sum_{k=0}^{\infty}
  \hbar^k c_k (a_1,a_2)
\]
satisfying the following properties:
\begin{enumerate}
\item
  Each one of the maps $c_k: A[[\hbar]]\otimes A[[\hbar]]\to A[[\hbar]]$ is
  $\complex[[\hbar]]$-bilinear;
\item
  One has $c_0(a_1, a_2)=a_1 \cdot a_2$ for all $a_1, a_2\in A$;
\item
  The relation
  $$
    a_1\star a_2 -c_0 (a_1,a_2) - \frac {i}{2} \hbar  \Pi (a_1, a_2)
    \in \hbar^2 A[[\hbar]]
  $$
  holds true for some representative $\Pi\in Z^2(A, A)$ of the Poisson
  structure and  all $a_1,a_2\in A$.
\end{enumerate}
\end{definition}

For an \'etale groupoid $\calg$ with an invariant symplectic form
$\omega$ and a invariant symplectic connection $\nabla$ on the
base, we define a Hochschild 2-cochain on $C^{\infty}(\calg)$ by
\begin{equation}
\label{def:pi}
  \Pi( a_1 , a_2)(g)= \sum_{g_1 \, g_2 =g} \,
  \pi(g )(da_1(g_1), da_2(g_2)), \quad g\in \calg \: \:
  a_1,a_2 \in C^{\infty}(\calg),
\end{equation}
where $da_1(g_1)$ and $da_2(g_2)$ have been pulled back to $g$
along the maps $t$ and $s$, and $\pi$ is the Poisson structure
associated to the symplectic form $\omega$. This definition is
legitimate because $t$ and $s$ are local diffeomorphisms. It was
proved \cite{t1:def-gpd} that this Hochschild 2-cochain gives rise
to a Poisson structure on $C^{\infty}(\calg)$ if there is an
invariant symplectic connection.

Tang \cite{t1:def-gpd} showed that the above noncommutative
Poisson structure $\Pi$ on the groupoid algebra admits a formal
deformation quantization. Such a deformation can be constructed as
follows: first using Fedosov's construction \cite{fe:book}, given
an invariant symplectic connection, we construct an invariant star
product on the algebra of smooth functions on the unit space
$\calg^{(0)}$. The deformation of the groupoid algebra
$C^{\infty}(\calg)$ is a crossed product algebra of the above
deformation on the base $C^{\infty}(\calg^{(0)})$ and the
associated pseudogroup $\calg$ action.

\subsection{Rankin-Cohen deformation}
It is well known that if $f(z)$ is a modular form, $\frac{1}{2\pi
i}\frac{d}{dz}f$ is not a modular form any more. Following
\cite{cm:modular-heck}, we introduce a differential operator $X$
as
$$
X\stackrel{def}{=}\frac{1}{2\pi i}\frac{d}{dz}-\frac{1}{12\pi
i}\frac{d}{dz}(\log\Delta)\cdot Y,
$$
where
$\Delta(z)=(2\pi)^{12}\eta^{24}(z)=(2\pi)^{12}q\prod_{n=1}^{\infty}(1-q^n)^{24},
\ q=e^{2\pi z}$ and $Y(f)=\frac{k}{2}f,\ \forall f\in \calm_k$,
the space of modular forms of weight $k$.

It is straightforward to check that $X$ and $Y$ acts on
$\calm=\oplus_k\calm_k$ satisfying $[Y, X]=X$. Under these two
operators, the Rankin-Cohen bracket $RC_n$ can be written as
follows, for $f\in \calm_k,\ g\in \calm_l$
$$
\begin{array}{ll}
RC_n(f, g)&=\sum_{r+s=n}(-1)^r\left(\begin{array}{c}n+k-1\\
s\end{array}\right)\left(\begin{array}{c}n+l-1\\
r\end{array}\right)f^{(r)}g^{(s)},
\end{array}
$$
where $f^{(r)}$ (or $g^{(s)})$ is the $r$-th (or $s$-th)
derivative of $f$ (or $g)$, and
$(\alpha)_k\stackrel{def}{=}\alpha(\alpha+1)\cdots (\alpha+k-1)$.

In \cite{z:deformation}, Zagier observed that $\sum_n RC_n$
defines an associative product on $\calm$. This product actually
defines a universal deformation formula of the Lie algebra $h_1$,
consisting of $X, Y$ with $[Y, X]=X$, since $h_1$ acts on $\calm$
injectively. It is worth mentioning that $h_1$ is the Lie algebra
of the $``ax+b"$ group.

Inspired by the Rankin-Cohen brackets, Connes and Moscovici
\cite{cm:deformation} introduced a family of Rankin-Cohen type
elements in $(\calh_1\otimes \calh_1)[[\hbar]]$ as follows.

\begin{dfn}
\label{dfn:proj} {\rm {(\cite{cm:deformation})}}Let $\calh_1$ act
on an algebra $A$. This action is called projective if
$\delta_2'\stackrel{def}{=}\delta_1^2-\frac{1}{2}\delta_2$ is
inner implemented by an element $\Omega \in A$, so that
$$
\delta_2'(a)=[\Omega, a],\hspace{1cm} \forall a\in A,
$$
and
$$
\delta_k(\Omega)=0, \hspace{1cm} \forall k \in \mathbb{N}.
$$
\end{dfn}
Assume that the action of $\calh_1$ action an algebra $A$ is
projective. Define
\begin{equation}
\label{eq:cm-def}
\begin{array}{ll}
RC&=\sum_{n=0}^{\infty}\hbar^{n}\sum_{k=0}^{n}\frac{A_k}{k!}(2Y+k)_{n-k}\otimes\frac{B_{n-k}}{(n-k)!}(2Y+n-k)_k\\
A_{m+1}&=S(X)A_m-m\Omega^0(Y-\frac{m-1}{2})A_{m-1},\\
B_{m+1}&=XB_m-m\Omega(Y-\frac{m-1}{2})B_{m-1},
\end{array}
\end{equation}
where $\Omega^0$ is the right multiplication of $\Omega$.

Connes and Moscovici \cite{cm:deformation} proved that $RC$
defines a universal deformation formula of a projective $\calh_1$
action.
\section{Universal deformation of $h_1$}
If we set all $\delta_n$ to be $0$, the Lie algebra $H_1$ is
reduced to $h_1$, the Lie algebra of the ``$ax+b$" group, and
$\calh_1$ becomes $\calu(h_1)$, the universal enveloping algebra
of $h_1$. In this case, $RC$ defined by (\ref{eq:cm-def}) is
simplified to the following universal deformation formula of
$h_1$,
\begin{equation}\label{RC}
RC_{n}(a,b)\stackrel{def}{=}\sum_{k=0}^n\left[
\frac{(-1)^k}{k!}X^k(2Y+k)_{n-k}(a)\;\frac{1}{(n-k)!}X^{n-k}(2Y+n-k)_{k}(b)
\right],
\end{equation}
where $X,Y\in h_1$ are such that $[Y,X]=X$,
$(\alpha)_k\stackrel{def}{=}\alpha(\alpha+1)...(\alpha+k-1)$, and
$a,b\in A$.

We spend this section studying this universal deformation.
\subsection{Giaquinto-Zhang's deformation of $h_1$}
A nice deformation formula for $h_1$ has already been given by
Giaquinto and Zhang  \cite{gzh:ax+b}[Thm 2.20]: Given two elements
$X,Y$ with $[Y,X]=X$, the following expression defines a universal
deformation formula(UDF) of the Hopf algebra associated to $h_1$
\[
F=\sum_{n=0}^{\infty}\frac{t^n}{n!}F_n=1\ts 1+t X\wedge Y
+\frac{t^2}{2!}\left( X^2\otimes Y_{2} -2 XY_{1}\otimes XY_{1}
+Y_{2}\otimes X^2 \right)+\cdots,
\]
where $F_n$ is defined to be $ F_n=\sum _{r=0}^n (-1)^r {n\choose
r}X^{n-r}Y_{r}\otimes X^{r}Y_{n-r}$.

\begin{prop}
\label{prop:gzh}The above defined $F$ can be realized by the
standard Moyal product.
\end{prop}

$\pf$ We consider the space $\mathbb R\times \mathbb R_+$ on which
$X$ and $Y$ act as $Y=-y \frac{\partial}{\partial y}$, and
$X=\frac{1}{y} \frac{\partial}{\partial x}$. It is obvious that the
action of $X$ and $Y$ on $\reals \times \reals_+$ is injective.

With the following identity,
$$
Y_{r}=Y(Y+1)\cdots(Y+r-1)=(-y)^r \frac{\partial^r}{\partial y^r},
$$
it is straightforward to check that the above defined $F$ in this
representation is equal to the Moyal product. $\Box$

\subsection{Rankin-Cohen deformation of $h_1$}
We should point out that the above universal deformation formula
of $h_1$ is not equal to the one induced from $RC$ in Equation
(\ref{RC}). However, we will show that it is equivalent to the
Giaquinto-Zhang's deformation.

We set $(V,\omega):=(\R^2=\{(p,q)\},dp\wedge dq)$ and denote by
$\h=\h(V,\omega):=V\times\R$ the associated Heisenberg algebra.
Setting $\frakg:=\mathfrak{sl}_2(\R)=\SPAN_\R\{H,E,F\}$,
$([H,E]=2E,\;[H,F]=-2F,\; [E,F]=H)$, we form the natural
semi-direct product $\tilde{\frakg}:=\frakg\times\h$. The
(infinitesimal) affine linear action $\tilde{\g}\to\Gamma(T(V))$
is then strongly hamiltonian. We let $\lambda:\tilde{\frakg}\to
C^\infty(V)$ denote the corresponding moment map. Explicitly,
denoting fundamental vector fields by $A^\star_x:=\ddto\exp(-t
A)\cdot x \quad A\in\tilde{\frakg}$, one has
\[
\begin{array}{lllll}
H^\star=-p\partial_p+q\partial_q;& E^\star=-q\partial_p;&
F^\star=-p\partial_q;& P^\star=-\partial_p;&
Q^\star=-\partial_q;\\
\lambda_H=pq;& \lambda_E=\frac{1}{2}q^2;&
\lambda_F=-\frac{1}{2}p^2;&\lambda_P=q;&\lambda_Q=-p.
\end{array}
\]
We have that $[A^\star,B^\star]=[A,B]^\star$ and
$\lambda_{[A,B]}=\{\lambda_A,\lambda_B\}$ where
$\{u,v\}=\partial_pu\partial_qv-\partial_pv
\partial_qu$, and $A, B\in \tilde{\frakg}$.

Let $S:=AN=\exp(\SPAN\{H,E\})$ denote the Iwasawa component in
$SL(2,\R)$, which is the $``ax+b"$ group. We consider the open
orbit $\CO\stackrel{def}{=}S\cdot(0,1)$ in $V$, which is equal to
the set $[q>0]$. Since $S$ acts simply transitively on $\CO$, we
have the identification $\phi:S\to\CO:g\mapsto g\cdot(1,0)$. We
still denote by $\lambda:\tilde{\frakg}\to C^{\infty}(S)$ the
transported restricted moment map, that is:
\begin{equation}
\lambda_A :=
\phi^\star(\lambda_A|_{\CO})\qquad(A\in\tilde{\frakg}).
\end{equation}
\begin{lem}\label{DIAG}
Denoting by $\tilde{X}_g:=\ddto \;g\exp(tX)$ the left-invariant
vector field associated to $X\in h_1=\mbox{Lie }(S)$, one has:
\begin{enumerate}
\item[(i)]
$\tilde{H}\;.\;\lambda_{X+v}\;=\;(-2)\;\lambda_X\;+\;(-1)\;
\lambda_v\quad\mbox{for }X\in\frakg\mbox{ and }v\in V;$
\item[(ii)]$\tilde{ E}^r\;.\;\lambda_X\;=\;0\mbox{ for }r\geq 3, $
for all $X\in \frakg$; \item[(iii)]
$\tilde{E}^r\;.\;\lambda_v\;=\;0\mbox{ for }r\geq 2, $ for all
$v\in V$.
\end{enumerate}
\end{lem}
\Pf A convenient parametrization of the group manifold $S$ is
given by:
$$
\R^2\to S:(a,\ell)\mapsto\exp(aH)\exp(\ell E).
$$
In these coordinates, the group law reads
$(a,\ell)\cdot(a',\ell')=(a+a',e^{-2a'}\ell+\ell')$. We deduce the
expressions for the left-invariant vector fields:
$$
\tilde{H}=\partial_a-2\ell\partial_\ell\;;\;\tilde{E}=\partial_\ell.
$$
The corresponding chart on the orbit $\CO\simeq S$ is given by
$$
p=e^a\ell\;;\;q=e^{-a}.
$$
Note that this is a global Darboux chart on $\CO$ as for $da\wedge
d\ell=\pm\phi^\star\omega|_{\CO}$. The corresponding (uncomplete)
moment map reads as
$$
\lambda_H=\ell\;;\;\lambda_E=\frac{1}{2}e^{-2a}\;;\;\lambda_F=-\frac{1}{2}\ell^2e^{2a}\;;\;
\lambda_P=e^{-a}\;;\;\lambda_Q=-e^a\ell.
$$
A straightforward computation then yields the lemma. $\Box$

From (\ref{RC}), for any left ${\cal U}(h_1)$ action on an algebra
$A$, the Rankin-Cohen brackets on $\calu(h_1)$ is defined by,
\[
RC_{n}(a,b):=\sum_{k=0}^n\left[
\frac{(-1)^k}{k!}X^k(2Y+k)_{n-k}(a)\;\frac{1}{k!}X^{n-k}(2Y+n-k)_{k}(b)
\right],
\]
where $X,Y\in h_1$ are such that $[Y,X]=X$,
$(\alpha)_k\stackrel{def}{=}\alpha(\alpha+1)...(\alpha+k-1)$, and
$a,b\in A$.

Since $h_1$ acts as left invariant vector fields on $S$, $\calu
(h_1)$ acts as left invariant differential operators on
$C^{\infty}(S)$, and $RC_n$, an element of $\calu (h_1)\otimes
\calu(h_1)$, acts as a left invariant bidifferential operator on
$C^{\infty}(S)$. Since $[H,E]=2E$, we set
\[
\tilde{H}=2Y\mbox{ and }\tilde{E}=X.
\]
\begin{lem}
\label{lem:moment-quant} For all $A$ in $\tilde{\frakg}$, we have
\begin{equation}
[\lambda_A\;,\;u]_n\stackrel{def}{=}
RC_n(\lambda_A\;,\;u)-RC_n(u\;,\;\lambda_A) \;=\;0\quad\mbox{for
}\quad n\neq1.
\end{equation}
\end{lem}
\Pf For $X\in\frakg$ and $v\in V$, Lemma \ref{DIAG} implies that
$X^k(2Y+r)_{s}.\lambda_{X+v}=(-2+r)_sX^k\lambda_X+(-1+r)_sX^k\lambda_v=0$
if $k>2$. Therefore, in the expression (\ref{RC}) of
$RC_n(\lambda_{X+v},u)$ only the first three terms corresponding
to $k=0,1,2$ contribute. In each of them the following (left hand
side) factor occurs:
\begin{equation}\label{I}
\bullet \mbox{ for }k=0: (-2)_n\lambda_X+(-1)_n\lambda_v\;;
\end{equation}
\begin{equation}\label{II}
\bullet \mbox{ for }k=1:
\tilde{E}.[(-1)_{n-1}\lambda_X+(0)_{n-1}\lambda_v]\;;
\end{equation}
\begin{equation}\label{III}
\bullet \mbox{ for }k=2:
\tilde{E}^2.[(0)_{n-2}\lambda_X+(1)_{n-2}\lambda_v].
\end{equation}
\begin{enumerate}
\item The first expression (\ref{I}) vanishes identically for
$n\geq3$. Indeed, $(-2)_{n}=(-2)(-2+1)(-2+2)...(-2+n-1)$ is zero
as soon as $n-1\geq2$; and similarly for $(-1)_n$;

\item In the same way, the second expression (\ref{II}) vanishes for
$n-2\geq1$, i.e. $n\geq3$;

\item At last, the third expression (\ref{III}) is equal to
$(n-2)!\tilde{E}^2(\lambda_v)$ which is identically zero by Lemma
\ref{DIAG} item  (iii). We conclude by observing that $RC_{0}$ and
$RC_{2}$ are symmetric. $\Box$
\end{enumerate}

By Lemma \ref{lem:moment-quant}, the Rankin-Cohen deformation
(\ref{eq:cm-def}) defines a $\tilde{\frakg}$ invariant star
product on $(V, \omega)$. In Corollary 2, Section 2.7 of
\cite{gutt:sl2}, Gutt showed that there is a unique
$\tilde{\frakg}$-invariant star product on $(V, \omega)$, which is
the standard Moyal product. We conclude that the Rankin-Cohen
deformation on $C^{\infty}(S)$ is identical to the Moyal product.

\begin{prop}
\label{prop:moyal} The reduced Rankin-Cohen deformation realized
on $\CO\subset V$ coincides with the restriction to $\CO$ of the
standard Moyal product on $(V,\Omega)$.
\end{prop}

To generalize the construction in Proposition \ref{prop:moyal}, we
explain its relation to Fedosov's construction of deformation
quantization of symplectic manifolds.

The natural action of $S\simeq``ax+b"$ on $\R$,
\[
\exp(aH+nE)\cdot x_{1}:=e^{2a}x_1+ne^a,
\]
lifts to $T^\star(\R)=\R^2$ as
\[
\exp(aH+nE)\cdot(x_{1},x_2):=(e^{2a}x_1+ne^a,e^{-2a}x_2).
\]
The $S$-orbit $\tilde{\CO}$ of point $\tilde{o}:=(0,1)=dx_1|_0\in
T^\star(\R^2)$ is then naturally isomorphic as $S$-homogeneous
space to $\CO\subset V$; namely one has the identification:
\[
\varphi:\CO\to\tilde{\CO}:g\cdot e_2\mapsto g\cdot \tilde{o}.
\]
In $(p,q)$-coordinates on $\CO$, this reads:
\[
\varphi(p,q)=\left(\;\frac{p}{2q}\;,\;q^2\;\right).
\]
Identifying $\tilde{\CO}$ with $S$ (via $\varphi\circ\phi$), we
obtain the expressions for the left invariant vector fields:
\[
\tilde{H}=-2x_2\partial_{x_2}\;;\;\tilde{E}=\frac{1}{x_2}\partial_{x_1}.
\]
In particular, we set
\[
\tilde{H}=2Y\mbox{ and }\tilde{E}=X.
\]
By letting $\nabla^\CO$ denote the restriction to $\CO$ of the
standard symmetric flat connection on $V$
($\nabla^\CO_{\partial_p}\partial_p=\nabla^\CO_{\partial_q}\partial_p=
\nabla^\CO_{\partial_q}\partial_q=0$), and setting
\[
\nabla^{\tilde{\CO}}:=\varphi(\nabla^\CO),
\]
we obtain a symplectic connection on $\tilde {\calo}$,
\begin{equation}
\label{eq:conn-flat}
\nabla^{\tilde{\CO}}_{\partial_{x_1}}\partial_{x_1}=0\;;\;\nabla^{\tilde{\CO}}_{\partial_{x_1}}\partial_{x_2}
=\frac{1}{2x_2}\partial_{x_1}\;;\;\nabla^{\tilde{\CO}}_{\partial_{x_2}}\partial_{x_2}=-
\frac{1}{2x_2}\partial_{x_2}.
\end{equation}

We identify $\tilde{\calo}$ with $\reals \times \reals^+$, and use
$\nabla^{\tilde{\calo}}$ to construct deformation quantization
 of $(\reals \times \reals ^+, \omega\stackrel{def}{=}dx\wedge dy)$ as described in Section 2.2.

\begin{cor}
\label{prop:star-fedosov}The reduced Rankin-Cohen deformation on
$\tilde{\calo}$ is identical to Fedosov's construction of the star
product on $(\tilde{\calo}, \omega)$ using the connection
$\nabla^{\tilde{\calo}}$ with the characteristic form equal to
$\frac{1}{i\hbar}\omega$.
\end{cor}

\section{Projective structures}
To reconstruct Connes-Moscovici's Rankin-Cohen deformation, we
need to understand the geometric meaning of their Definition
\ref{dfn:proj}, a projective structure.

\subsection{The flat case}
We look at the connection $\nabla^{\tilde{\calo}}$
 considered in Section 3, (\ref{eq:conn-flat}).

\begin{prop}
\label{prop:flat-conn} The connection $\nabla^{\tilde{\calo}}$
(\ref{eq:conn-flat}) is invariant under the local diffeomorphism
$\phi: x_1\mapsto \tilde{x_1}\stackrel{def}{=}\phi(x_1),
x_2\mapsto \tilde{x_2}\stackrel{def}{=} \frac{x_2}{\phi'(x_1)}$ if
and only if $\delta_2 '(\phi)=0$. Here $\calh_1$ acts on $\phi$ as
in Section 2.1.
\end{prop}
{\bf Notation:}We use $\nabla$ to replace $\nabla^{\tilde{\calo}}$
in the rest of the paper.\\
\noindent{$\pf$}We have the following transformation rules of
vector fields.
\[
\begin{array}{l}
\frac{\partial}{\partial
\tilde{x}_1}=\frac{1}{\phi'(x_1)}\frac{\partial}{\partial
x_1}+\frac{\phi''}{\phi'^2}x_2\frac{\partial}{\partial x_2},\\
\frac{\partial}{\partial
\tilde{x}_2}=\phi'\frac{\partial}{\partial x_2}.
\end{array}
\]

The invariance of $\nabla$ implies that we should have
\[
\begin{array}{ll}
\nabla_{\phi_*(\frac{\partial}{\partial
x_1})}\phi_*(\frac{\partial}{\partial x_1})&=
\nabla_{\phi'(x_1)\frac{\partial}{\partial
\tilde{x}_1}-\frac{\phi''}{\phi'^2}x_2\frac{\partial}{\partial
\tilde{x}_2}}({\phi'(x_1)\frac{\partial}{\partial
\tilde{x}_1}-\frac{\phi''}{\phi'^2}x_2\frac{\partial}{\partial
\tilde{x}_2}})\\
&=\phi'^2 \nabla_\frac{\partial}{\partial
\tilde{x}_1} \frac{\partial}{\partial \tilde{x}_1}+
\phi'\frac{\partial}{\partial
\tilde{x}_1}(\phi')\frac{\partial}{\partial
\tilde{x}_1}-\frac{\phi''}{\phi'}x_2\nabla_\frac{\partial}{\partial
\tilde{x}_1}\frac{\partial}{\partial
\tilde{x}_2}-\phi'\frac{\partial}{\partial
\tilde{x}_1}(\frac{\phi''}{\phi'^2}x_2)\frac{\partial}{\partial
\tilde{x}_2}\\
&-\frac{\phi''}{\phi'}x_2\nabla_\frac{\partial}{\partial
\tilde{x}_2}\frac{\partial}{\partial
\tilde{x}_1}-\frac{\phi''}{\phi'^2}x_2\frac{\partial}{\partial
\tilde{x}_2}(\phi')\frac{\partial}{\partial
\tilde{x}_1}+(\frac{\phi''}{\phi'^2}x_2)^2\nabla_\frac{\partial}{\partial
\tilde{x}_2}\frac{\partial}{\partial
\tilde{x}_2}\\
&+\frac{\phi''}{\phi'^2}x_2\frac{\partial}{\partial
\tilde{x}_2}(\frac{\phi''}{\phi'^2}x_2)\frac{\partial}{\partial \tilde{x}_2}\\
&=\phi'\frac{1}{\phi'}(\phi'')\frac{\partial}{\partial
\tilde{x}_1}-\phi'\frac{\phi''}{\phi'^2}x_2\frac{1}{2\tilde{x}_2}\frac{\partial}{\partial
\tilde{x}_1}
-\phi'[\frac{1}{\phi'}\frac{\phi'''\phi'^2-2\phi''^2\phi'}{(\phi'^2)^2}x_2+(\frac{\phi''}{\phi'^2})^2x_2]
\frac{\partial}{\partial \tilde{x}_2}\\
&
-\phi'\frac{\phi''}{\phi'^2}x_2\frac{1}{2\tilde{x}_2}\frac{\partial}{\partial
\tilde{x}_1}+0+(\frac{\phi''}{\phi'^2}x_2)^2\frac{1}{2\tilde{x}_2}\frac{\partial}{\partial
\tilde{x}_2}+\frac{\phi''}{\phi'^2}x_2 \phi'
\frac{\phi''}{\phi'^2}\frac{\partial}{\partial
\tilde{x}_2}\\
&=
-\frac{\phi'''\phi'-\frac{3}{2}\phi''^2}{\phi'^3}x_2\frac{\partial}{\partial
\tilde{x}_2},
\end{array}\\
\]
\[
\begin{array}{l}
\begin{array}{ll}
\nabla_{\phi_*(\frac{\partial}{\partial
x_1})}\phi_*(\frac{\partial}{\partial x_2}) & =
\nabla_{\phi'(x_1)\frac{\partial}{\partial
\tilde{x}_1}-\frac{\phi''}{\phi'^2}x_2\frac{\partial}{\partial
\tilde{x}_2}} (\frac{1}{\phi'}\frac{\partial}{\partial
\tilde{x}_2})\\
& =\phi' \frac{1}{\phi'}\nabla_\frac{\partial}{\partial
\tilde{x}_1} \frac{\partial}{\partial
\tilde{x}_2}+\phi'\frac{\partial}{\partial
\tilde{x}_1}(\frac{1}{\phi'})\frac{\partial}{\partial
\tilde{x}_2}-\frac{\phi''}{\phi'^2}x_2
\frac{1}{\phi'}\nabla_{\frac{\partial}{\partial \tilde{x}_2}}
\frac{\partial}{\partial
\tilde{x}_2}-\frac{\phi''}{\phi'^2}x_2\frac{\partial}{\partial
\tilde{x}_2}(\frac{1}{\phi'})\frac{\partial}{\partial
\tilde{x}_2}\\
& =\frac{1}{2\tilde{x}_2}\frac{\partial}{\partial
\tilde{x}_1}+\phi'\frac{1}{\phi'}(-\frac{\phi''}{\phi'^2})\frac{\partial}{\partial
\tilde{x}_2}-\frac{\phi''}{\phi'^2}x_2
\frac{1}{\phi'}(-\frac{1}{2\tilde{x}_2}\frac{\partial}{\partial \tilde{x}_2})-0\\
& =\frac{1}{2\tilde{x}_2}\frac{\partial}{\partial
\tilde{x}_1}-\frac{1}{2}\frac{\phi''}{\phi'^2}\frac{\partial}{\partial
\tilde{x}_2}=\phi_ * (\frac{1}{2x_2}\frac{\partial}{\partial
x_1}),
\end{array}\\
\\
\begin{array}{ll}
\nabla_{\phi_*(\frac{\partial}{\partial
x_2})}\phi_*(\frac{\partial}{\partial x_2})  &
\nabla_{\frac{1}{\phi'}\frac{\partial}{\partial \tilde{x}_2}}
(\frac{1}{\phi'}\frac{\partial}{\partial
\tilde{x}_2})=\frac{1}{\phi'^2}\nabla_{\frac{\partial}{\partial
\tilde{x}_2}} (\frac{\partial}{\partial
\tilde{x}_2})+\frac{1}{\phi'}\frac{\partial}{\partial
\tilde{x}_2}(\frac{1}{\phi'})\frac{\partial}{\partial
\tilde{x}_2}\\
& =
\frac{1}{\phi'^2}(-\frac{1}{2\tilde{x}_2})\frac{\partial}{\partial
\tilde{x}_2}+0=\phi_ * (-\frac{1}{2x_2}\frac{\partial}{\partial
x_2}){\Big{|}}_{(\tilde{x}_1,\tilde{x}_2)}.
\end{array}
\end{array}
\]
We see easily that the invariance of the connection under $\phi$
is equivalent to $\phi'''\phi'-\frac{3}{2}\phi''^2=0$, i.e.
$\delta_2'(\phi)=0$.\hspace{1cm} $\Box$

\subsection{The general case}
For the general case of nontrivial $\delta_2'$, we look at the
following connection.
\begin{equation}
\label{eq:conn-gen}
\begin{array}{ll}
\nabla_{\frac{\partial}{\partial x_1}}\frac{\partial}{\partial
x_1}=\mu(x_1, x_2) \frac{\partial}{\partial x_2} , &
\nabla_{\frac{\partial}{\partial x_1}}\frac{\partial}{\partial
x_2}=\frac{1}{2x_2}\frac{\partial}{\partial x_1},\\
& \\
\nabla_{\frac{\partial}{\partial x_2}}\frac{\partial}{\partial
x_1}=\frac{1}{2x_2}\frac{\partial}{\partial x_1}, &
\nabla_{\frac{\partial}{\partial x_2}}\frac{\partial}{\partial
x_2}=-\frac{1}{2x_2}\frac{\partial}{\partial x_2}.
\end{array}
\end{equation}
Here $\mu$ is a suitable function.

\begin{thm}
\label{thm:proj}Let $\Gamma$ be a pseudogroup generated by local
diffeomorphisms on $\reals$ acting on $\reals\times \reals^+$ by
$\phi: x_1\mapsto \phi(x_1), x_2\mapsto \frac{x_2}{\phi'(x_1)}$,
$\forall \phi\in \Gamma$. Assume that the dimension of the fixed
point set of each element $\phi\in \Gamma$ is strictly less than
2. The connection $\nabla$ in (\ref{eq:conn-gen}) is invariant
under $\Gamma$ if and only if the $\calh_1$ action on the
corresponding groupoid algebra $\Gamma\ltimes
C_c^{\infty}(\reals\times \reals^+)$ is projective.
\end{thm}
\medskip
$\pf$ \noindent Given a local diffeomorphism $\phi$, we have the
following quantity different from the proof of Proposition
\ref{prop:flat-conn}. All the others are same.

$$
\begin{array}{ll}
\nabla_{\phi_*(\frac{\partial}{\partial
x_1})}\phi_*(\frac{\partial}{\partial x_1}) & =
\nabla_{\phi'(x_1)\frac{\partial}{\partial
\tilde{x}_1}-\frac{\phi''}{\phi'^2}x_2\frac{\partial}{\partial
\tilde{x}_2}}({\phi'(x_1)\frac{\partial}{\partial
\tilde{x}_1}-\frac{\phi''}{\phi'^2}x_2\frac{\partial}{\partial
\tilde{x}_2}})\\
& =  \phi'^2 \nabla_\frac{\partial}{\partial \tilde{x}_1}
\frac{\partial}{\partial \tilde{x}_1}+
\phi'\frac{\partial}{\partial
\tilde{x}_1}(\phi')\frac{\partial}{\partial
\tilde{x}_1}-\frac{\phi''}{\phi'}x_2\nabla_\frac{\partial}{\partial
\tilde{x}_1}\frac{\partial}{\partial
\tilde{x}_2}\\
&-\phi'\frac{\partial}{\partial
\tilde{x}_1}(\frac{\phi''}{\phi'^2}x_2)\frac{\partial}{\partial
\tilde{x}_2}
-\frac{\phi''}{\phi'}x_2\nabla_\frac{\partial}{\partial
\tilde{x}_2}\frac{\partial}{\partial
\tilde{x}_1}-\frac{\phi''}{\phi'^2}x_2\frac{\partial}{\partial
\tilde{x}_2}(\phi')\frac{\partial}{\partial
\tilde{x}_1}\\
&+(\frac{\phi''}{\phi'^2}x_2)^2\nabla_\frac{\partial}{\partial
\tilde{x}_2}\frac{\partial}{\partial
\tilde{x}_2}+\frac{\phi''}{\phi'^2}x_2\frac{\partial}{\partial
\tilde{x}_2}(\frac{\phi''}{\phi'^2}x_2)\frac{\partial}{\partial
\tilde{x}_2}\\
&= \phi'^2\mu(\tilde{x}_1, \tilde{x}_2)\frac{\partial}{\partial
\tilde{x}_2}+\phi'\frac{1}{\phi'}(\phi'')\frac{\partial}{\partial
\tilde{x}_1}-\phi'\frac{\phi''}{\phi'^2}x_2\frac{1}{2\tilde{x}_2}\frac{\partial}{\partial
\tilde{x}_1}\\
&-\phi'[\frac{1}{\phi'}\frac{\phi'''\phi'^2-2\phi''^2\phi'}{(\phi'^2)^2}x_2+(\frac{\phi''}{\phi'^2})^2x_2]
\frac{\partial}{\partial \tilde{x}_2}
-\phi'\frac{\phi''}{\phi'^2}x_2\frac{1}{2\tilde{x}_2}\frac{\partial}{\partial
\tilde{x}_1}\\
&+(\frac{\phi''}{\phi'^2}x_2)^2\frac{1}{2\tilde{x}_2}\frac{\partial}{\partial
\tilde{x}_2}+\frac{\phi''}{\phi'^2}x_2 \phi'
\frac{\phi''}{\phi'^2}\frac{\partial}{\partial
\tilde{x}_2}\\
&= [\phi'^2\mu(\tilde{x}_1,
\tilde{x}_2)-\frac{\phi'''\phi'-\frac{3}{2}\phi''^2}{\phi'^3}x_2]\frac{\partial}{\partial
\tilde{x}_2}. \end{array}
$$

By the invariance of $\nabla$, we have
\[
[\phi'^2\mu(\tilde{x}_1,
\tilde{x}_2)-\frac{\phi'''\phi'-\frac{3}{2}\phi''^2}{\phi'^3}x_2]\frac{\partial}{\partial
\tilde{x}_2}=\phi_*(\mu(x_1)\frac{\partial}{\partial
x_2})=\mu(x_1, x_2)\frac{1}{\phi'}\frac{\partial}{\partial
\tilde{x}_2},
\]
and
\begin{equation}
\label{eq:delta2'}
\frac{\phi'''\phi'-\frac{3}{2}\phi''^2}{\phi'^3}x_2=\phi'^2\mu(\phi(x_1),
\frac{x_2}{\phi'})-\frac{1}{\phi'}\mu(x_1, x_2).
\end{equation}

By Equation (\ref{eq:delta2'}), we have
\begin{equation}
\label{eq:inner}
\frac{\phi'''\phi'-\frac{3}{2}\phi''^2}{\phi'^2}x_2^2=\phi'^4\tilde{x}_2\mu(\phi(x_1),
\frac{x_2}{\phi'})-x_2\mu(x_1, x_2).
\end{equation}
\begin{enumerate}
\item $\Rightarrow$. Let $\phi$ be an element in $\Gamma$.\\
We introduce $\nu=\frac {\mu(x_1, x_2)}{x_2}$, and Equation
(\ref{eq:inner}) is equivalent to
\[
\frac{\phi'''\phi'-\frac{3}{2}\phi''^2}{\phi'^2}=\phi'^2
\nu(\phi(x_1), \frac{x_2}{\phi'})-\nu(x_1, x_2).
\]

Define $\omega(x_1,x_2)=\nu(x_1, \frac{1}{x_2})$, and we have
\[\frac{\phi'''\phi'-\frac{3}{2}\phi''^2}{\phi'^2}\\
=\phi'^2 \nu(\phi(x_1), \frac{x_2}{\phi'})-\nu(x_1, x_2)\\
=\phi'^2\omega(\phi(x_1), \frac{\phi'}{x_2})-\omega(x_1,
\frac{1}{x_2}).\]

Introduce $y=\frac{1}{x_2}$, the above equation gives
\begin{equation}
\label{eq:omega}
\frac{\phi'''\phi'-\frac{3}{2}\phi''^2}{\phi'^2}=\phi'^2\omega(\phi(x_1),
\phi' y)-\omega(x, y).
\end{equation}

Finally, letting $\Omega(x,y)=y^2\omega(x,y),\ x_1=x$, we see that
Equation (\ref{eq:omega}) implies
\[
\frac{\phi'''\phi'-\frac{3}{2}\phi''^2}{\phi'^2}y^2=\phi'^2y^2\omega(\phi(x_1),
\phi' y)-\omega(x, y)y^2=(\phi^{-1})^*(\Omega)(x,y)-\Omega(x,y).
\]

The left hand side of the above equation is equal to the
expression of $\delta_2'(\phi^{-1})$. The above equality shows
that $\delta_2'$ is inner when we consider the $\calh_1$ action on
the foliation groupoid $FX\rtimes \calg$ as in Section 2.1.
\item $\Leftarrow$. Suppose that the $\calh_1$ action on
$\Gamma\ltimes C_c^{\infty}(\reals\times \reals^+)$ is projective.\\
We first show that if the $\calh_1$ action is projection on
$\Gamma \ltimes C_c^{\infty}(\reals\times \reals^+)$, the support
of $\Omega$ has to be on the unit space. We write
$\Omega=\sum_{\alpha\in \Gamma}\Omega_{\alpha}U_{\alpha}$ and
$\delta_2'(U_\phi)U_{\phi}=[\Omega, U_\phi]$, and have the
following observations.
\begin{enumerate}
\item From $\delta_i(\Omega)=0,\ \forall i>0$, we know that
$\delta_{i}(U_{\alpha})\Omega_{\alpha}=0, \forall \alpha$.
\item
From $\delta_i(f)=0$ for any $f\in C_c^{\infty}(\reals\times
\reals^+)$, we have that $[\Omega, f]=\sum_{\alpha\in
\Gamma}(\alpha^{*}(f)-f)\Omega_{\alpha} U_{\alpha}$. Therefore
$(\alpha^{*}(f)-f)\Omega_{\alpha}=0$, for all $\alpha\in \Gamma$.
\end{enumerate}

For a given $\alpha\in \Gamma$ not equal to identity, we have that
$\delta_{i}(U_{\alpha})\Omega_{\alpha}=0,\ \forall i>0$ and
$(\alpha^{*}(f)-f)\Omega_{\alpha}=0$. If there is $x_0\in
\reals\times \reals^+ $ such that $\Omega_{\alpha}(x_0)\ne 0$,
then at $x_0$, there is a neighborhood $N$ of $x_0$ on which
$\delta_i(U_\alpha)=0$. In particular
$\delta_1(U_{\alpha})=\log((\alpha^{-1})^{'})'=0$. Solving this
differential equation, we know that $\alpha$ on $N$ must act like
$\alpha:(x_1, x_2)\mapsto (ax_1+b, ax_2 )$. By the fact that
$(\alpha^*(f)-f)\Omega_{\alpha}(x_0)=0$ on $N$, for any smooth
function, we know that $\alpha(x_0)=x_0$. The same argument show
that all $x\in N$ has to be fixed by $\alpha$, since
$\Omega_{\alpha}(x)\ne0$. But this contradicts our assumption that
the fixed point set of $\alpha$ is at most 1 dimensional. This
shows that $\Omega_{\alpha}=0$.

From the above argument, we know that $\Omega$ has to be supported
on the unit space. At this time, the projective condition is
equivalent to
\[
\delta_2(\phi^{-1})=y^2\frac{\phi'''\phi'-\frac{3}{2}{\phi''}^2}{\phi'^2}U_{\phi}=(\Omega-\phi^*(\Omega))U_{\phi}.
\]

From (\ref{eq:omega}) and the transformation there, we know that
the existence of $\Omega$ implies the existence of an invariant
connection like (\ref{eq:conn-gen}). $\Box$
\end{enumerate}
\begin{rmk}
\label{rmk:frame}Here, for calculation convenience, we have
identified the Frame bundle $F\reals$ with the cotangent bundle
$T^* \reals$ by $\tau: (x,y)\mapsto (x, \frac{1}{y})$. The
connection $\nabla$ is defined on $T^* \reals$. By $\tau$, it is
also defined on $F\reals$.
\end{rmk}

In Theorem \ref{thm:proj}, the assumption that the fixed point set
of any element in $\Gamma$ is at most one dimensional is only used
in the sufficient part of the proof. Generally, $\Omega$ is
supported on the fixed point set $B^{(0)}$ of $\Gamma$, i.e.
$\{(\gamma, x)|\ \gamma\in \Gamma, \gamma(x)=x\}$. $\Gamma$ acts
on $B^{(0)}$, by conjugation action. The similar result of Theorem
\ref{thm:proj} is extended to this general situation without any
extra effort.
\\
\noindent{{\bf Theorem 4.2'} }{\em Let $\Gamma$ be a pseudogroup
generated by local diffeomorphisms on $\reals$ and $B^{(0)}=\{
(\gamma, x)\in \Gamma\times \reals\times \reals^+| \gamma\cdot
x=x\}$ be the fixed point set. The projective action $(\rho,
\Omega)$ of $\calh_1$ on $\Gamma\ltimes C_c^{\infty}(\reals\times
\reals^+)$ is one to one correspondent to a $\Gamma$ invariant
connection $\nabla$ on $\reals\times \reals^+$ of form
(\ref{eq:conn-gen}) and a smooth function $f$ on $\Gamma\times
\reals\times \reals^+$, which is supported on $B^{(0)}-\{(id,
x)|x\in \reals\times \reals^+\}$ and invariant under $\Gamma$
conjugation action.}
\section{Universal deformation formula for $\calh_1$}
In this section, we will use a Fedosov type construction to
reconstruct the universal deformation formula of $\calh_1$
originally constructed by Connes and Moscovici
\cite{cm:deformation}.
\subsection{Zagier's deformation}

In this subsection, we discuss the influence of the above new
connection (\ref{eq:conn-gen}) on the star product
(\ref{eq:star}).

\begin{cor}
\label{cor:curvature} The connection $\nabla$ (\ref{eq:conn-gen})
is flat if and only if $\mu(x_1,x_2)=x_2\nu(x_1)$, where
$\nu(x_1)$ is an arbitrary smooth function on $\reals$.
\end{cor}
$\pf$ The curvature of $\nabla$ can be directly calculated to be
equal to
\[
\begin{array}{ll}
R(\frac{\partial}{\partial x_1}, \frac{\partial}{\partial
x_2})(\frac{\partial}{\partial
x_1})&=(\frac{\mu}{x_2}-\frac{\partial \mu}{\partial
x_2})\frac{\partial}{\partial x_2}\\
R(\frac{\partial}{\partial x_1}, \frac{\partial }{\partial
x_2})(\frac{\partial}{\partial x_2})&=0.
\end{array}
\]
Therefore, $R=0$ if and only if $\frac{\mu}{x_2}-\frac{\partial
\mu}{\partial x_2}=0$. The solution of this first order
differential equation is that $\mu=x_2\nu(x_1)$, where $\nu(x_1)$
is an arbitrary smooth function on $\reals$. $\Box$

In this section, we restrict ourselves to the case that the
connection (\ref{eq:conn-gen}) is flat, which means that $\mu(x_1,
x_2)=x_2\nu(x_1)$. We consider the deformation quantization of
$(\reals\times \reals^+, dx_1\wedge dx_2)$ using this connection.

The Christoffel symbols of the connection $\nabla^{\tilde{\calo}}$
are calculated as follows,
\[
\Gamma_{11}^1=\Gamma_{12}^2=\Gamma_{21}^2=\Gamma_{22}^1=0,
\Gamma_{11}^2=\mu, \Gamma_{12}^1=\Gamma_{21}^1=\frac{1}{2x_2},
\Gamma_{22}^2=-\frac{1}{2x_2}.
\]

Taking the Formula (5.1.8) in \cite{fe:book} with the same
notations, we have

$$\begin{array}{ll}
\Gamma_{111}=\omega_{11}\Gamma_{11}^1+\omega_{12}\Gamma_{11}^2=\omega_{12}\mu,&
\Gamma_{211}=\omega_{21}\Gamma_{11}^1+\omega_{22}\Gamma_{11}^2=0,\\
\Gamma_{112}=\omega_{11}\Gamma_{12}^1+\omega_{12}\Gamma_{12}^2=0,&
\Gamma_{121}=\omega_{11}\Gamma_{21}^1+\omega_{12}\Gamma_{21}^2=0,\\
\Gamma_{212}=\omega_{21}\Gamma_{12}^1+\omega_{22}\Gamma_{12}^2=\frac{1}{2x_2}\omega_{21},&
\Gamma_{221}=\omega_{21}\Gamma_{21}^1+\omega_{22}\Gamma_{21}^2=\frac{1}{2x_2}\omega_{21},\\
\Gamma_{122}=\omega_{11}\Gamma_{22}^1+\omega_{12}\Gamma_{22}^2=-\frac{1}{2x_2}\omega_{12},&
\Gamma_{222}=\omega_{21}\Gamma_{22}^1+\omega_{22}\Gamma_{22}^2=0.
\end{array}$$

We have the following expression for $\Gamma,\ \Gamma\circ a,\
a\circ \Gamma$, and $[\Gamma, a]$.
$$
\Gamma=\frac{1}{2}\omega_{21}\{[-\mu(u^1)^2+\frac{1}{2}(2
u^2)^2]dx_1+\frac{1}{2} 2 u^1 u^2 dx_2\},
$$
and
$$
\begin{array}{ll}
\frac{i}{h}[\Gamma, a]&=\sum (\frac{1}{2}(-\mu)2a_{m,n}(u^1)^m n
(u^2)^{n-1}-\frac{1}{4x_2}2 a_{m,n}m(u^1)^{m-1} (u^2)^{n+1})dx_1\\
&+\frac{1}{4x_2}(2 a_{m,n}(u^1)^m n (u^2)^n-2 a_{m,n}m(u^1)^m
(u^2)^n)d x_2.
\end{array}
$$

It is a direct check that when $\mu=x_1\nu(x_2)$, $\nabla^2$ and
$D^2$ are both 0. By Theorem \ref{thm:fedosov}, for each $f\in
C^{\infty}(\reals\times \reals_+)[[\hbar]]$, there is a unique
solution of the equation $Da=0$ with $a_{0,0}=f$. In the
following, we calculate the explicit expression of $a$.

The expression of $Da$ is calculated as follows.
$$
\begin{array}{ll}
D a&=\partial a - \delta a= -\delta a + d a +
\frac{i}{h}[\Gamma,a]\\
& =-\sum a_{m,n}m (u^1)^{m-1} (u^2)^n d x_1 -\sum a_{m,n} (u^1)^m
n (u^2)^{n-1} d x_2\\
&+\sum \frac{\partial a_{m,n}}{\partial x_1}(u^1)^m (u^2)^n d
x_1+\sum \frac{\partial a_{m,n}}{\partial x_2} (u^1)^m (u^2)^n d
x_2\\
&+[-\mu\sum a_{m,n} n (u^1)^{m+1} (u^2)^{n-1}-\sum \frac{a_{m,n}}{2x_2}m(u^1)^{m-1} (u^2)^{n+1}] d x_1\\
&+ \sum \frac{a_{m,n}}{2x_2}(n-m)(u^1)^m (u^2)^n d x_2.
\end{array}
$$

The equation $Da=0$ gives the following system of differential
equations:
$$- a_{m+1,n}(m+1)+ \frac{\partial a_{m,n}}{\partial
x_1}-(n+1)\mu a_{m-1,n+1}-\frac{a_{m+1,n-1}}{2x_2}(m+1)=0,$$ and
$$-a_{m,n+1}(n+1)+ \frac{\partial a_{m,n}}{\partial
x_2}+\frac{a_{m,n}}{2x_2}(n-m)=0.$$

Given $a_{0,0}=f$, we solve the system of equations by induction.

$$
\begin{array}{ll}
a_{m,0}&=\frac{1}{m}(\frac{\partial a_{m-1,0}}{\partial x_1}-\mu
a_{m-2,1})=\frac{1}{m}(\frac{\partial a_{m-1,0}}{\partial x_1}-\mu
(\frac{\partial}{\partial x_2
}-\frac{m-2}{2x_2})a_{m-2,0}),\\
a_{m,n}&=\frac{1}{n!}(\frac{\partial}{\partial x_2
}-\frac{m}{2x_2})\cdots(\frac{\partial}{\partial x_2
}+\frac{n-m-1}{2x_2})a_{m,0}.
\end{array}
$$

If we set
$$
\begin{array}{ll}
X&=\frac{1}{x_2}\frac{\partial}{\partial x_1},\\
Y&=-x_2\frac{\partial}{\partial x_2},\\
\end{array}
$$
it is direct check that
$$
\begin{array}{ll}
A_{m+1}&=-X A_m -
m\frac{\mu}{x_2^3}(Y-\frac{m-1}{2})A_{m-1},\\
B_{m+1}&=X B_m -
m\frac{\mu}{x_2^3}(Y-\frac{m-1}{2})B_{m-1},\\
a_{m,n}&=\frac{(-1)^nx_2^{m-n}}{n!}\frac{A_m}{m!}(Y+\frac{m}{2})\cdots(Y+\frac{m+n-1}{2})a,\\
b_{n,m}&=\frac{(-1)^mx_2^{n-m}}{m!}\frac{B_n}{n!}(Y+\frac{n}{2})\cdots(Y+\frac{m+n-1}{2})b.
\end{array}
$$
The above expression of $A_m,B_m$ is exactly identical to the
recurrence relation as described in (2.9) of [2] of Connes and
Moscovici with $S(X)=-X$, and
$\Omega=\frac{\mu}{x^3_2}=\frac{\nu}{x_2^2}$. The star product
constructed in this way defines the Zagier's deformation
\cite{z:deformation} for $h_1$ constructed from Rankin-Cohen
brackets on modular forms with a forth degree element.
\begin{rmk}
For computation reasons, we have chosen that a special form of
connections defined by Equation (\ref{eq:conn-gen}), which is
flat. Because of the flatness, the calculation is quite simple and
transparent. When the connection is not flat, Fedosov's
construction still works, but the calculation is much more
complicated. However, the star product should be able to be
expressed by the same formula.
\end{rmk}
\begin{rmk}
As explained in Remark \ref{rmk:frame}, the connection and the
star product discussed in this subsection are both on the
cotangent bundle $T^* \reals$. However, all these constructions
can be pulled back to the frame bundle by $\tau$(See Remark
\ref{rmk:frame}) without any difficulty.
\end{rmk}

\subsection{Full injectivity}
We have shown in the last subsection that the deformation
quantization of the standard symplectic structure on the upper
half plane using the connection (\ref{eq:conn-gen}) with $\mu(x_1,
x_2)=x_2\nu(x_1)$ gives rise to Zagier's deformation formula on
modular forms. To generalize this deformation to a universal
deformation formula of a projective $\calh_1$ action, we adapt the
method used by Connes and Moscovici \cite{cm:deformation}[Sec. 3]
to our situation. We briefly recall their construction in the
following, and refer to \cite{cm:deformation} for the detail.

Firstly, we introduce a free abelian algebra $P$ with a set of
generators indexed by $\integers_{\geq 0}$, $Z_0, Z_1, \cdots,
Z_n, \cdots$. On $P$, we define a $\calh_1$ action as follows,
\[
Y(Z_j)\stackrel{def}{=}(j+2)Z_j, X(Z_j)\stackrel{def}{=}Z_{j+1},
\delta_k(p)=0,\  \forall p\in P,\ j\geq0.
\]

Secondly, we consider the crossed product algebra
$\tilde{\calh_1}\stackrel{def}{=}P\rtimes \calh_1\ltimes P$, which
is equal to $P\otimes \calh_1\otimes P$ as a vector space. Denote
this algebra  by $\tilde{\calh_1}$. Connes and Moscovici defines
on $\tilde{\calh_1}$ an Hopf algebra structure over $P$, with
$\alpha, \beta: P\to \tilde{\calh_1}$ defined by
\[
\alpha(p)=p\rtimes 1\ltimes 1,\ \ \ \ \ \beta(q)=1\rtimes 1\ltimes
q,\ \ \forall p, q \in P.
\]

Thirdly, to deal with the projective structure,  we define
$\tilde{\delta_2'}\stackrel{def}{=}\delta_2-\frac{1}{2}\delta_2-\alpha(Z_0)+\beta(Z_0)$,
$\tilde{\calh_s}$ as the quotient of $\tilde{\calh_1}$ by the
ideal generated by $\tilde{\delta_2}'$. $\tilde{\calh}_s$ is still
a Hopf algebra over $P$ because
$\Delta(\tilde{\delta_2'})=\tilde{\delta_2}'\otimes 1+1\otimes
\tilde{\delta_2}'$.

Fixing a function $\mu(x_1, x_2)$, we consider a pseudogroup
$\Gamma$ action on $\reals$ whose lifting onto $T^*\reals$
preserves the connection $\nabla$ (\ref{eq:conn-gen}) defined by
$\mu$. By Theorem \ref{thm:proj}, the $\calh_1$ action on the
corresponding groupoid algebra $\cala_{ \mu,
\Gamma}\stackrel{def}{=}C^{\infty}_c(\reals\times
\reals^+)\rtimes\Gamma$ is projective with $\Omega$ defined in the
proof.

We define $\rho_{\mu, \Gamma}:P\to \cala_{\mu, \Gamma}$ by
$\rho(Z_k)=X^k(\Omega)$ and make $\cala_{\mu, \Gamma}$ into a
module algebra over  $\tilde{\calh}_1|P$ by
\[
\chi_{\mu, \Gamma}(p\rtimes h\rtimes
q)(U_{\gamma}f)\stackrel{def}{=}\rho_{\mu,
\Gamma}(p)h(U_{\gamma}f)\rho_{\mu, \Gamma}(q).
\]

One easily checks that $\cala_{\mu, \Gamma}$ becomes a module
algebra over $\tilde{\calh}_s|P$ because when the $\calh_1$ action
is projective, $\tilde{\delta'_2}$ acts as $0$.

We define action $\chi_{\mu,\Gamma}^n$,
$$
\chi_{\mu,
\Gamma}^{(n)}:\underbrace{\tilde{\calh}_s\otimes_{P}\cdots \otimes
\tilde{\calh}_s} _n\to
\call(\underbrace{\cala_{\mu,\Gamma}\otimes\cdots\cala_{\mu,
\Gamma}}_{n}, \cala_{\mu, \Gamma})
$$
by means of acting on each components, where $\call$ means the set
of linear maps.

We fix $\mu=x_1\nu(x_1)$, and have the following Proposition
analogous to \cite{cm:deformation}[Prop. 12].
\begin{prop}
\label{prop:injective} For each $n\in \mathbb{N}$,
$\bigcap_{\nu(x_1), \Gamma }Ker\chi_{x_2\nu(x_1),
\Gamma}^{(n)}=0$.
\end{prop}
$\pf$ There is no difference between the proofs for different $n$.
Therefore, for simplicity, we only prove the proposition for
$n=1$.

Following the proof of \cite{cm:deformation}[Prop. 12], an
arbitrary element of $\tilde{\calh}_s$ can be written uniquely as
a finite sum of the form
\[
H=\sum_{j,k,l,m}\alpha(p_{jklm})\beta(q_{jklm})\delta_1^jX^kY^l,
\]
where $p,q\in P$.

Let $\chi_{x_2\nu(x_1), \Gamma}(H)=0$, for arbitrary $\nu(x_1)$
and pseudogroup $\Gamma$ preserving the connection defined by
$x_2\nu(x_1)$. From the proof of Theorem \ref{thm:proj}, we know
that in this case, $\Omega=x_2^2\nu(x_1)$.

If $U_{\gamma}f\in \cala_{x_2\nu(x_2), \Gamma}$, then
$$
\sum_{j,k,l,m}\rho_{x_1\nu(x_2),
\Gamma}(p_{jklm})\gamma^*(\rho_{x_1\nu(x_2)}(q_{jklm}))\delta_1(\gamma)^jX^kY^l(f).
$$

We notice that $f$ can be arbitrary smooth function on
$\reals\times\reals^+$, and
$X^kY^l=x_2^{m+l}\frac{d^k}{dx_1^m}\frac{d^l}{dx_2^l}$. This
implies that
\[
\sum_{j,m}\rho_{x_1\nu(x_2),
\Gamma}(p_{jklm})\gamma^*(\rho_{x_1\nu(x_2)}(q_{jklm}))\delta_1(\gamma)^j=0,
\]
for any $l,m$.

To prove the Proposition, we consider the following family of
algebras,  $\cala_{x_2\nu(x_2), \Gamma}$.

Fix a diffeomorphism $\phi_{O_1,O_2}$ from an open set $O_1\subset
\reals$ to the other open set $O_2\subset \reals$, with $O_1$
disjoint from $O_2$. The disjointness between $O_1$ and $O_2$
makes the set $\Gamma_{\phi}\stackrel{def}{=}\{id|_{\reals},
id|_{O_1}, id|_{O_2}, \phi, \phi^{-1},\}$ into a pseudogroup.
Starting with any connection $\nabla_1$ of the form
(\ref{eq:conn-gen}) with $\mu=x_2\nu(x_1)$ on $O_1$, we first push
forward this connection to $O_2$ by $\phi$, and then extend the
connections defined on $O_1$ and $O_2$ to a global connection
$\tilde{\nabla}$ on $\reals\times \reals^+$. The extension of the
connection is well defined because $O_1$ is disjoint from $O_2$,
(we may need to restrict to a smaller open subset $O_2'$ of $O_2$
by a cutoff function) and is $\Gamma_{\phi}$ invariant by its
definition. According to our construction, we have that
$\tilde{\calh_s}$ act on the corresponding groupoid algebra
$\cala_{\phi_{O_1, O_2}, \tilde{\nabla}}$.

Now at any $x\in \reals$, we fix $O_1$ containing $x$, and let
$O_2$, $\phi$, $\nabla_1$ vary. It is not hard to see that if $H$
vanishes on this family of algebra $\cala_{\phi_{O_1, O_2},
\tilde{\nabla}}$, we must have that $H$ vanishes at $x$, because
$H$ has only finite number of terms but this family of algebras
has infinitely many freedoms. Hence $H$ has to be equal to 0.
$\Box$

\subsection{Universal deformation $\calh_1$ with a projective structure}
We consider the groupoid algebra $\cala_{x_2\nu(x_1), \Gamma}$.
Because the connection defined by $x_1\nu(x_1)$ in
(\ref{eq:conn-gen}) is $\Gamma$ invariant, the results in Section
2.3 implies that the symplectic form $\frac{{\rm d}x\wedge{\rm
d}y}{y^2}$ on $\reals\times \reals_+$, which is invariant under
any $\Gamma$, defines a noncommutative Poisson structure on
$C_c^{\infty}(\reals\times \reals^+)\rtimes \Gamma$. Furthermore,
we extend this Poisson structure to a deformation of
$C_c^{\infty}(\reals\times \reals^+)\rtimes \Gamma$. This
deformation can be realized by the crossed product of the star
product constructed in Section 5.1 with $\Gamma$.

In Section 5.1, the $\star$ product is expressed as follows: for
$f, g\in C_c^{\infty}(\reals\times \reals_+)$,
$$
\begin{array}{ll}
f\star
g&=\sum_{n=0}^{\infty}\hbar^{n}\sum_{k=0}^{n}\frac{A_k}{k!}(2Y-k)_{n-k}(a)\frac{B_{n-k}}{(n-k)!}(2Y-n+k)_k(b)\\
A_{m+1}&=-X A_m -
mx_2\mu(Y-\frac{m-1}{2})A_{m-1}=-XA_m-m\Omega (Y-\frac{m-1}{2})A_{m-1},\\
B_{m+1}&=X B_m -
mx_2\mu(Y-\frac{m-1}{2})B_{m-1}=XB_m-m\Omega(Y-\frac{m-1}{2})B_{m-1}.
\end{array}
$$
The crossed product of $\star$ with $\Gamma$ is written as
$f_{\gamma}U_\gamma \ast
g_{\beta}U_\beta\stackrel{def}{=}f_{\gamma}\star
\gamma^*(g_{\beta})U_{\gamma\beta}$ defines a deformation
quantization of $ C^{\infty}_c(\reals\times \reals^+)\rtimes
\Gamma$.

According to the formulas of $\star$ and the $\Gamma$ crossed
product, the deformed product $\ast$ on $
C^{\infty}_c(\reals\times \reals^+)\rtimes \Gamma$ can be
expressed by $\tilde{\calh}_s$ as follows,
\[
\begin{array}{ll}
RC&=\sum_{n=0}^{\infty}\hbar^{n}\sum_{k=0}^{n}\frac{A_k}{k!}(2Y+k)_{n-k}\otimes\frac{B_{n-k}}{(n-k)!}(2Y+n-k)_k\\
A_{m+1}&=S(X)A_m-m\Omega^0(Y-\frac{m-1}{2})A_{m-1},\\
B_{m+1}&=XB_m-m\Omega(Y-\frac{m-1}{2})B_{m-1},
\end{array}
\]
where $\Omega^0$ is the right multiplication of $\Omega$.

By Proposition \ref{prop:injective}, we conclude $RC$ can be
pulled back to $\tilde{\calh}_s$ and defines an associative
universal deformation for any projective $\calh_1$ actions.

\section{Deformation without Projective structures---noncommutative Poisson structure}
In the above deformation (\ref{eq:cm-def}), we have assumed the
action to be projective. One can ask whether one can go beyond
this. Recently, a construction of Bressler, Gorokhovsky, Nest, and
Tsygan strongly suggests that this general RC deformation may
still exist.

In this section, we look at the first order approximation of the
general deformation. We prove that $RC_1$ generally defines a
noncommutative Poisson structure without any assumptions.

\begin{prop}
\label{prop:ncpoi}For an $\calh_1$ action on an $A$,
$RC_1=-X\otimes 2Y+2Y\otimes X+\delta_1Y\otimes 2Y$ defines a
noncommutative Poisson structure on $A$.
\end{prop}

$\pf$ The proof of this proposition is calculation. We need to
find an element $B$ in $\calh_1\otimes \calh_1$, such that for any
$a, b, c\in A$,
\[
aB(b,c)-B(ab, c)+B(a, bc)-B(a,b)c=RC_1(RC_1(a, b), c)-RC_1(a,
RC_1(b, c)).
\]

In order to find such a $B$, we first look at the special case
where the Hopf algebra action is projective. In this case, the
associativity of the Connes-Moscovici's universal deformation
formula of $\calh_1$ implies that $RC_2$ is a right choice of $B$.

For a general $\calh_1$ action, we first look at the following
term
$$B'=S(X)^2\otimes Y(2Y+1)+S(X)(2Y+1)\otimes X(2Y+1)+Y(2Y+1)\otimes X^2.$$
We calculate the difference between the Hochschild coboundary of
$B'$ and $[RC_1, RC_1]$.

$$\begin{array}{ll}
&(b(B')-[RC_1, RC_1])(a,b,c)\\
=& 4 Y a \delta_2'Y b Y c+2Y^2a\delta_2'b
Yc+2Ya\delta_2'bYc+2Ya\delta_2'bY^2c\\
\\
=&-2 [a\delta_2'Y^2 b Y c-\delta_2' Y^2(ab) Yc+\delta_2'Y^2 a
Y(bc)-\delta_2' Y^2 a (Yb) c]-4\delta_2'YaYbYc\\
&-2\delta_2'aY^2b Yc+2Ya\delta_2'aY^2c+2Ya\delta'_2aYc\\
\\
=& -2 [a\delta_2'Y^2 b Y c-\delta_2' Y^2(ab) Yc+\delta_2'Y^2 a
Y(bc)-\delta_2' Y^2 a (Yb) c]\\
&-2[a\delta_2'YbY^2c-\delta_2'Y(ab)Y^2c+\delta_2'YaY^2(bc)-\delta_2'Ya(Y^2b)c]\\
&-2\delta_2'aY^2bYc-2\delta_2'aYbY^2c+2Ya\delta_2'bYc\\
\\
=& -2 [a\delta_2'Y^2 b Y c-\delta_2' Y^2(ab) Yc+\delta_2'Y^2 a
Y(bc)-\delta_2' Y^2 a (Yb) c]\\
&-2[a\delta_2'YbY^2c-\delta_2'Y(ab)Y^2c+\delta_2'YaY^2(bc)-\delta_2'Ya(Y^2b)c]\\
&-\frac{2}{3}[a\delta_2'bY^3c-\delta_2'(ab)Y^3c+\delta_2'aY^3(bc)+\delta_2'aY^3(b)c]+2Ya\delta_2'bYc\\
=& -2 [a\delta_2'Y^2 b Y c-\delta_2' Y^2(ab) Yc+\delta_2'Y^2 a
Y(bc)-\delta_2' Y^2 a (Yb) c]\\
&-2[a\delta_2'YbY^2c-\delta_2'Y(ab)Y^2c+\delta_2'YaY^2(bc)-\delta_2'Ya(Y^2b)c]\\
&-\frac{2}{3}[a\delta_2'bY^3c-\delta_2'(ab)Y^3c+\delta_2'aY^3(bc)+\delta_2'aY^3(b)c]\\
&-2[a\delta_2'YbYc-\delta_2'Y(ab)Yc+\delta_2'YaY(bc)-\delta_2'YaY(b)c]\\
&-[a\delta_2'bY^2c-\delta_2'(ab)Y^2c+\delta_2'aY^2(bc)-\delta_2'aY^2(b)c],
\end{array}$$
where $b(B')$ is the Hochschild coboundary of $B'$ and
$\delta_2'=\delta_2-\frac{1}{2}\delta_1^2$.

It is straightforward to check the following identities.
$$
\begin{array}{rl}
b(\delta_2'Y^2\otimes Y)(a, b, c)&=a\delta_2'Y^2 b Y c-\delta_2'
Y^2(ab) Yc+\delta_2'Y^2 a Y(bc)-\delta_2' Y^2 a (Yb) c,\\
b(\delta_2'\otimes Y^3)(a, b, c)&=a\delta_2' b Y^3 c-\delta_2'(ab)
Y^3 c+\delta_2' a Y^3(bc)-\delta_2' a (Y^3 b) c,\\
b(\delta_2'Y\otimes Y)(a, b, c)&=a \delta_2'Y b Y c-\delta_2'
Y(ab)Yc+\delta_2' Y aY(bc)-\delta_2' Y a(Yb)c,\\
b(\delta_2'\otimes Y^2)(a, b, c)&=a\delta_2' b Y^2 c
-\delta_2'(ab)Y^2 c +\delta_2' a Y^2(bc)-\delta_2' a Y^2 b
c\\
b(\delta_2' Y\otimes
Y^2)(a,b,c)&=a\delta_2'YbY^2c-\delta_2'Y(ab)Y^2c+\delta_2'YaY^2(bc)-\delta_2'YaY^2(b)c.
\end{array}
$$

Therefore, the calculation suggests the introduction of
$B''=2\delta_2'+Y^2\otimes Y+\frac{2}{3}\delta_2'\otimes
Y^3+2\delta_2'Y\otimes Y^2+2\delta_2'Y\otimes Y+\delta_2'\otimes
Y^2$ and $B=B'+B''$. And we have $b(B)=b(B'+B'')=[RC_1, RC_1]$.
$\Box$

\noindent{Pierre Bieliavsky}, D\'epartement de Math\'ematique,
Universit\'e Catholique de Louvain, Chemin du cyclotron, 2, 1348
Louvain-La-Neuve, Belgium. Email: bieliavsky@math.ucl.ac.be.
\vspace{2mm}

\noindent{Xiang Tang}, Department of Mathematics, Washington University, St. Louis, MO, 63130, U.S.A.,  Email:xtang@math.wustl.edu.

\vspace{2mm}
\noindent{Yijun Yao}, Centre de Math\'ematiques
\'Ecole Polytechnique, Ecole Polytechnique, 91128 Palaiseau Cedex,
France, Email: yao@math.polytechnique.fr.
\end{document}